\documentclass[reqno,10pt]{amsart}

\usepackage{amssymb}
\usepackage{amsmath}
\usepackage{color}
\usepackage{graphicx,epsfig}
\usepackage{mathrsfs}
\usepackage{amstext}
\usepackage{latexsym}

\usepackage[labelformat=empty]{caption}
\usepackage{enumerate}

\usepackage[pdftex]{hyperref}
\usepackage[pdftex,open]{bookmark}
\theoremstyle{plain}
\newtheorem{theorem}{Theorem} %[section]
\newtheorem{corollary}[theorem]{Corollary}

\newtheorem{example}[theorem]{Example}

\newtheorem{proposition}[theorem]{Proposition}
\theoremstyle{remark}
\newtheorem{remark}{Remark}
\newcommand{\reel}{\mathbb{R}}
\newcommand{\nat}{\mathbb{N}}

\newcommand{\ds}{\displaystyle}

\newcommand{\abs}[1]{\left\vert #1\right\vert }

\newcommand{\bg}{\medskip\goodbreak}

\begin{document}
\title[Elementary evaluation of  $\ds\int_0^\infty\frac{\sin^pt}{t^q}d t$]
{Elementary evaluation of  $\ds\int_0^\infty\frac{\sin^pt}{t^q}d t$}
\author[Omran Kouba]{Omran Kouba$^\dag$}
\address{Department of Mathematics \\
Higher Institute for Applied Sciences and Technology\\
P.O. Box 31983, Damascus, Syria.}
\email{\href{mailto:omran_kouba@hiast.edu.sy}{omran\_kouba@hiast.edu.sy}}
\thanks{$^\dag$ Department of Mathematics, Higher Institute for Applied Sciences and Technology.}

%%%%%%%%%%%%%%%%%%%%%%%%%%%%
\begin{abstract}
Let $p$ and $q$ be two positive integers. The goal of this note is to demonstrate, in a very simple and
elementary way and without using advanced tools, a formula that  expresses the value of the integral 
$\int_0^\infty\frac{\sin^pt}{t^q}dt$ when it converges.\par

\end{abstract}

\maketitle
%%%%%%%%%%%%%%%%%%%%%%%%
\section{\sc Introduction }\label{sec1}

For positive integers $p$ and $q$ we consider the family of integrals
\begin{equation*}
I(p,q)=\int_0^\infty\frac{\sin^pt}{t^q}dt
\end{equation*}
It seems that these integrals were first considered by N. I. Lobachevskii \cite{lob}.
An explicit evaluation of $I(p,q)$ when $q-p$ is even was given by T. Hayashi \cite{hay}. The formula given in \cite{hay} shows that in this case $I(p,q)$ is a rational multiple of $\pi$ and this was precisely the object of Problem 11423 proposed to the American Mathematical Monthly \cite{mon}. A detailed evaluation of these integrals can be found in the literature, for example an evaluation using distribution theory and Fourier transforms can be found in \cite{kol}. Other methods using contour integration can also be applied to derive these formulas. 

The aim of this note is to present elementary proofs for the formulas for $I(p,q)$ when the integral converges.

\section{\sc The main results}\label{sec2}

Our basic tool in the proof is the  sequence of polynomials $(P_n)_{n\geq 0}$ defined as follows:
\begin{align*}
	P_{2n}(X)=&\sum_{k=0}^{n}\frac{(-1)^k}{(2k)!}X^{2k}\\
	P_{2n+1}(X)=&\sum_{k=0}^{n}\frac{(-1)^k}{(2k+1)!}X^{2k+1}
	\end{align*}
and the associated sequence of functions $(f_n)_{n\geq1}$ given by:
\begin{align*}\forall n\geq0,\qquad &f_{2n+1}:\reel\longrightarrow\reel, f_{2n+1}(t)=(-1)^n(\sin t-P_{2n-1}(t))\\
	&f_{2n+2}:\reel\longrightarrow\reel, f_{2n+2}(t)=(-1)^{n+1}(\cos t-P_{2n}(t))
	\end{align*}
 with the convention $P_{-1}=0$.

The evaluation of the integral $I(p,q)$ is based upon the next proposition which summarizes some properties of the functions  $(f_n)_{n\geq1}$:\bg

\begin{proposition} \label{pr1} The sequence of functions  $(f_n)_{n\geq1}$ satisfies the following properties:

\begin{enumerate} [$i$.]
\item For every $n\in\nat^*$ we have $f^\prime_{n+1}=f_n$.
\item For every $n\in\nat^*$ we have $\ds\lim_{t\to0}\frac{f_{n}(t)}{ t^n}=\frac{1}{n!}$.
\item For every $n\in\nat^*$ the improper integral  $\ds J_n=\int_0^\infty \frac{f_n(t)}{ t^n}d t$ converges, and \[J_n=\frac{1}{(n-1)!}J_1=\frac{1}{(n-1)!}\cdot    \frac{\pi}{2}\]
\item  For every $n\geq2$ the improper integral  $\ds K_n=\int_1^\infty\frac{f_{n-1}(t)}{t^n}d t$    converges, and for all $\lambda>0$ we have:
\[\lim_{X\to\infty}\int_X^{\lambda X}\frac{f_{n+1}(t)}{t^n}d    t=\frac{\ln\lambda}{(n-1)!}\quad\text{and}\quad
\lim_{X\to\infty}\int_1^{X}\frac{f_{n+1}(t)}{t^n}d t=+\infty\]
\item  For every  $(m,q)\in\nat^2$ such that $1\leq q\leq m$ we have 
\begin{equation*}\label{Eq1}
\forall t\in\reel,\qquad(\sin t)^{2m}=\frac{1}{ 2^{2m-1}}\sum_{k=1}^m\binom{
2m}{ m-k}(-1)^{k+q}f_{2q}(2kt)\tag{1}
\end{equation*}
\item  For every  $(m,q)\in\nat^2$ such that $0\leq q\leq m$ we have 
\begin{equation*}
\label{Eq2}
\forall t\in\reel,\qquad(\sin t)^{2m+1}=\frac{1}{ 2^{2m}}\sum_{k=0}^m\binom{2m+1}{                 m-k}(-1)^{k+q}f_{2q+1}((2k+1)t)\tag{2}
\end{equation*}
\end{enumerate}
\end{proposition}

\begin{proof}

 $i.$ The verification of the first property is straightforward.\par
$ii.$ The second property follows from the power series expansions:
$$\cos t=P_{2n}(t)+\sum_{k=n+1}^\infty\frac{(-1)^k}{(2k)!}t^{2k}\quad\hbox{and}\quad
\sin t=P_{2n-1}(t)+\sum_{k=n}^\infty\frac{(-1)^k}{(2k+1)!}t^{2k+1}$$
which are valid for all $n\in\nat$ and all $t\in\reel$.\par
$iii.$  It is well-known that the integral $\ds J_1=\int_0^\infty\frac{\sin t}{ t}d t$ converges and that $\ds J_1=\frac{\pi}{2}$. In the case $n\geq2$ we have $\ds\frac{\abs{f_n(t)}}{ t^n}=O\left(\frac{1}{ t^2}\right)$ and this, with $ii$, proves the convergence of the integral $J_n$. Now, for $n\geq2$ and $X>0$ we have
\begin{align*}\int_0^X\frac{f_n(t)}{ t^n}d t=&\left[-\frac{f_n(t)}{ (n-1)t^{n-1}}\right]_0^X+\frac{1}{ n-1}
	\int_0^X\frac{f^\prime_n(t)}{ t^{n-1}}d t\\
	=&\frac{1}{ n-1}
	\int_0^X\frac{f_{n-1}(t)}{ t^{n-1}}d t+O\left(\frac{1}{ X}\right)
	\end{align*}
Letting $X$ tend to infinity we find that $\ds J_n=\frac{1}{n-1}J_{n-1}$ and this proves $iii$ by induction.\par

$iv.$ Considering two cases according to the parity of $n$ we see immediately that:
\begin{equation*}
\label{Eq3}
\forall t>0,\qquad\frac{f_{n+1}(t)}{ t^n}=\frac{1}{ (n-1)!}\cdot\frac{1}{ t}-\frac{f_{n-1}(t)}{ t^n}\tag{3}
\end{equation*}

But, for $n\geq2$,  we have $f_{n-1}(t)=O\left(t^{\max(n-3,0)}\right)$ in the neighborhood of $+\infty$, so  $\ds\frac{\abs{f_{n-1}(t)}}{t^n}=O\left(\frac{1}{  t^2}\right)$, and this proves the convergence of
$\ds K_n=\int_1^\infty\frac{f_{n-1}(t)}{ t^n}d t$. Now, using \eqref{Eq3} we conclude that
$\ds\lim_{X\to\infty}\int_1^{X}\frac{f_{n+1}(t)}{ t^n}d t=+\infty$ and that
$$\int_X^{\lambda X}\frac{f_{n+1}(t)}{t^n}=\frac{\ln \lambda}{(n-1)!}-\int_X^{\lambda X}\frac{f_{n-1}(t)}{ t^n}$$
The convergence of $K_n$ proves that $\ds\lim_{X\to\infty}\int_X^{\lambda X}\frac{f_{n-1}(t)}{ t^n}=0$ so that 
$$\lim_{X\to\infty}\int_X^{\lambda X}\frac{f_{n+1}(t)}{t^n}=\frac{\ln \lambda}{ (n-1)!}$$\par
$v.$ Let us start with a well-known and standard calculation. From Euler's formula : $\ds \sin t=\frac{e^{i t} - e^{-i t}}{ 2i}$ and using the binomial theorem we can write the following:
\begin{align*}
(\sin t)^{2m}=&\frac{(-1)^m}{ 2^{2m}}\left(\sum_{k=0}^{2m}\binom{2m}{k}(-1)^ke^{i kt}e^{-i(2m-k)t}\right)=
	\frac{(-1)^m}{ 2^{2m}}\left(\sum_{k=0}^{2m}\binom{2m}{k}(-1)^ke^{2i(k-m)t}\right)\\
	=&\frac{(-1)^m}{ 2^{2m}}\left(\sum_{k=0}^{m-1}\binom{2m}{k}(-1)^ke^{2i(k-m)t}+(-1)^m\binom{2m}{ m}  
	+\sum_{k=m+1}^{2m}\binom{2m}{k}(-1)^ke^{2i(k-m)t}\right)\\
	=&\frac{(-1)^m}{ 2^{2m}}\left(\sum_{k=0}^{m-1}\binom{2m}{k}(-1)^ke^{2i(k-m)t}+(-1)^m\binom{2m}{ m}  
	+\sum_{k=0}^{m-1}\binom{2m}{k}(-1)^ke^{-2i(k-m)t}\right)\\
	=&\frac{1}{ 2^{2m}}\left(\binom{2m}{ m}+2\sum_{k=0}^{m-1}\binom{2m}{k}(-1)^{m-k}\cos({2(k-m)t})\right)\\
	=&\frac{\binom{2m}{ m}}{ 2^{2m}}+\frac{1}{2^{2m-1}}\sum_{k=1}^{m}\binom{2m}{m-k}(-1)^{k}\cos({2kt})
\end{align*}
Using that $\cos u=(-1)^q f_{2q}(u)+P_{2q-2}(u)$ for all $u\in\reel$, we obtain
$$\forall t\in\reel,\qquad(\sin t)^{2m}=Q_{m,q}^0(t)+\frac{1}{2^{2m-1}}\sum_{k=1}^{m}
\binom{2m}{m-k}(-1)^{k+q}f_{2q}({2kt})$$
with $\ds Q_{m,q}^0(t)=\frac{\binom{2m}{ m}}{ 2^{2m}}+
\frac{1}{2^{2m-1}}\sum_{k=1}^{m}\left({2m\atop m-k}\right)(-1)^{k}P_{2q-2}({2kt})$.
\par
Since, in a neighborhood of  $0$, we have,  $\forall k\in\{1,2,\ldots,m\},\, f_{2q}(2kt)=O(t^{2q})$ and $(\sin t)^{2m}=O(t^{2m})=O(t^{2q})$ (this results from the assumption $q\leq m$), we conclude that the polynomial $Q_{m,q}^0$, which is of degree at most $2q-2$, satisfies : $Q_{m,q}^0(t)=O(t^{2q})$ in a neighborhood of  $0$. This proves  that  $Q_{m,q}^0=0$, and establishes \eqref{Eq1}.\par
$vi.$ In a similar way, starting from Euler's formula and using the binomial theorem we can write :
\begin{align*}
(\sin t)^{2m+1}=&\frac{(-1)^{m+1}}{ 2^{2m+1}i}\left(\sum_{k=0}^{2m+1}\binom{2m+1}{ k}(-1)^ke^{i kt}e^{-i(2m+1-k)t}\right)\\
	=&
	\frac{(-1)^{m+1}}{ 2^{2m+1}i}\left(\sum_{k=0}^{2m+1}\binom{2m+1}{ k}(-1)^ke^{i(2(k-m)-1)t}\right)\\
	=&\frac{(-1)^{m+1}}{ 2^{2m+1}i}\left(\sum_{k=0}^{m}\binom{2m+1}{ k}(-1)^ke^{i(2(k-m)-1)t}+\sum_{k=m+1}^{2m+1}\binom{2m+1}{ k}(-1)^ke^{i(2(k-m)-1)t}\right)\\
	=&\frac{(-1)^{m+1}}{ 2^{2m+1}i}\left(\sum_{k=0}^{m}\binom{2m+1}{ k}(-1)^ke^{i(2(k-m)-1)t}-\sum_{k=0}^{m}\left({2m+1\atop k}\right)(-1)^ke^{i(2(m-k)+1)t}\right)\\
	=&\frac{(-1)^{m}}{ 2^{2m}}\sum_{k=0}^{m}\binom{2m+1}{ k}(-1)^{k}\sin({(2(m-k)+1)t})\\
	=&\frac{1}{2^{2m}}\sum_{k=0}^{m}\binom{2m+1}{m- k}(-1)^{k}\sin({(2k+1)t})
\end{align*}
But $\sin u=(-1)^q f_{2q+1}(u)+P_{2q-1}(u)$ for all
$u\in\reel,$, hence
$$\forall t\in\reel,\qquad(\sin t)^{2m+1}=Q_{m,q}^1(t)+\frac{1}{2^{2m}}\sum_{k=0}^{m}\binom{2m+1}{m- k}(-1)^{k+q}f_{2q+1}({(2k+1)t})$$
with $\ds Q_{m,q}^1(t)=
\frac{1}{2^{2m}}\sum_{k=0}^{m}\binom{2m+1}{m- k}(-1)^{k}P_{2q-1}({(2k+1)t})$.
\par
Since, in a neighborhood of  $0$, we have,  $\forall k\in\{1,\ldots,m\},\, f_{2q+1}((2k+1)t)=O(t^{2q+1})$ and $(\sin t)^{2m+1}=O(t^{2m+1})=O(t^{2q+1})$ (this results from the assumption $q\leq m$), we conclude that the polynomial $Q_{m,q}^1$, which is of degree at most $2q-1$, satisfies : $Q_{m,q}^1(t)=O(t^{2q+1})$ in a neighborhood of  $0$. This proves that  $Q_{m,q}^1=0$, and establishes \eqref{Eq2}.
\end{proof}

\bg

\begin{theorem} Consider $(q,m)\in\nat$.
\begin{enumerate}[$i$.]
\item  If $1\leq q \leq m$ then $\ds\int_0^\infty\frac{\sin^{2m}t}{ t^{2q}}d t=
\frac{\pi}{2^{2m}}\sum_{k=1}^m(-1)^{k+q}\binom{2m}{ m-k} \frac{(2k)^{2q-1}}{(2q-1)!}$.

\item   If $0\leq q \leq m$ then $\ds\int_0^\infty\frac{\sin^{2m+1}t}{ t^{2q+1}}d t=
\frac{\pi}{2^{2m+1}}\sum_{k=0}^m(-1)^{k+q}\binom{2m+1}{ m-k} \frac{(2k+1)^{2q}}{(2q)!}$.
\item  If $2\leq q \leq m$ then $\ds\int_0^\infty\frac{\sin^{2m}t}{ t^{2q-1}}d t=
\frac{1}{2^{2m-1}}\sum_{k=1}^m(-1)^{k+q}\binom{2m}{ m-k}\frac{(2k)^{2q-2}}{(2q-2)!}\ln k$.
\item   If $1\leq q \leq m$ then $\ds\int_0^\infty\frac{\sin^{2m+1}t}{ t^{2q}}d t=
\frac{1}{2^{2m}}\sum_{k=0}^m(-1)^{k+q}\binom{2m+1}{ m-k}\frac{(2k+1)^{2q-1}}{(2q-1)!}\ln(2k+1)$.
\end{enumerate}
\end{theorem}
\bg

\begin{proof} $i.$ Using formula \eqref{Eq1} from Proposition \ref{pr1}, we conclude that
\begin{align*}
\int_0^\infty\frac{\sin^{2m}t}{ t^{2q}}d t=&\frac{1}{ 2^{2m-1}}\sum_{k=1}^m(-1)^{k+q}\binom{2m}{ m-k}\int_0^\infty{f_{2q}(2kt)\over t^{2q}}d t\\
	=&\frac{1}{ 2^{2m-1}}\sum_{k=1}^m(-1)^{k+q}\binom{2m}{ m-k} (2k)^{2q-1}\int_0^\infty\frac{f_{2q}(u)}{ u^{2q}}d u\\
	=&\frac{1}{2^{2m-1}}\sum_{k=1}^m(-1)^{k+q}\binom{2m}{ m-k}(2k)^{2q-1}J_{2q}
	\end{align*}
Now, using $iii$ from Proposition \ref{pr1}, we conclude that
$$\int_0^\infty\frac{\sin^{2m}t}{ t^{2q}}d t
=\frac{\pi}{ 2^{2m}}\sum_{k=1}^m(-1)^{k+q}\binom{2m}{ m-k} \frac{(2k)^{2q-1}}{(2q-1)!}$$
\bg

$ii.$ Using formula \eqref{Eq2} from Proposition \ref{pr1}, we conclude that
\begin{align*}
\int_0^\infty\frac{\sin^{2m+1}t}{ t^{2q+1}}d t=&\frac{1}{ 2^{2m}}\sum_{k=0}^m(-1)^{k+q}\binom{2m+1}{ m-k} \int_0^\infty\frac{f_{2q+1}((2k+1)t)}{ t^{2q+1}}d t\\
	=&\frac{1}{ 2^{2m}}\sum_{k=0}^m(-1)^{k+q}\binom{2m+1}{ m-k}(2k+1)^{2q}\int_0^\infty\frac{f_{2q+1}(u)}{ u^{2q+1}}d u\\
	=&\frac{1}{ 2^{2m}}\sum_{k=0}^m(-1)^{k+q}\binom{2m+1}{ m-k} (2k+1)^{2q}J_{2q+1}
\end{align*}
Now, using $iii$ from Proposition \ref{pr1}, we conclude that once more that
$$\int_0^\infty\frac{\sin^{2m+1}t}{ t^{2q+1}}d t
=\frac{\pi}{ 2^{2m+1}}\sum_{k=0}^m(-1)^{k+q}\binom{2m+1}{ m-k}\frac{(2k+1)^{2q}}{(2q)!}$$

\bg
$iii.$ Using formula \eqref{Eq1} from Proposition \ref{pr1}, we conclude that, for $X>0$ we have:
\begin{align*}\int_0^X\frac{\sin^{2m}t}{ t^{2q-1}}d t=&\frac{1}{ 2^{2m-1}}\sum_{k=1}^m(-1)^{k+q}\binom{2m}{ m-k}\int_0^X\frac{f_{2q}(2kt)}{ t^{2q-1}}d t\\
	=&\frac{1}{ 2^{2m-1}}\sum_{k=1}^m(-1)^{k+q}\binom{2m}{ m-k}(2k)^{2q-2}\int_0^{2kX}\frac{f_{2q}(u)}{ u^{2q-1}}d u\\
	=&\frac{1}{ 2^{2m-1}}\sum_{k=1}^m(-1)^{k+q}\binom{2m}{ m-k}(2k)^{2q-2}\int_{2X}^{2kX}\frac{f_{2q}(u)}{ u^{2q-1}}d u+C_{m,q}^0\int_{0}^{2X}\frac{f_{2q}(u)}{ u^{2q-1}}d u
	\end{align*}
with $\ds C_{m,q}^0=\frac{1}{ 2^{2m-1}}\sum_{k=1}^m(-1)^{k+q}\left({2m\atop m-k}\right)(2k)^{2q-2}$.
But, for $2\leq q\leq m$ the integral $\ds \int_0^\infty\frac{\sin^{2m}t}{ t^{2q-1}}d t$ converges, so that using $iv $ from
Proposition \ref{pr1}, we conclude that we must have $ C_{m,q}^0=0$, and consequently
$$\int_0^\infty\frac{\sin^{2m}t}{ t^{2q-1}}d t
=\frac{1}{ 2^{2m-1}}\sum_{k=1}^m(-1)^{k+q}\binom{2m}{ m-k} \frac{(2k)^{2q-2}}{(2q-2)!}\ln k$$
\bg

$iv.$ Using \eqref{Eq2} we conclude that, for $X>0$, we have :
\begin{align*}
	\int_0^X\frac{\sin^{2m+1}t}{ t^{2q}}d t=&\frac{1}{ 2^{2m}}\sum_{k=0}^m(-1)^{k+q}\binom{2m+1}{ m-k} \int_0^X\frac{f_{2q+1}((2k+1)t)}{ t^{2q}}d t\\
	=&\frac{1}{ 2^{2m}}\sum_{k=0}^m(-1)^{k+q}\binom{2m+1}{ m-k} (2k+1)^{2q-1}\int_0^{(2k+1)X}\frac{f_{2q+1}(u)}{ u^{2q}}d u\\
	=&\frac{1}{ 2^{2m}}\sum_{k=0}^m(-1)^{k+q}\binom{2m+1}{ m-k}(2k+1)^{2q-1}\int_{X}^{(2k+1)X}\frac{f_{2q+1}(u)}{ u^{2q}}d u+C_{m,q}^1\int_{0}^{X}\frac{f_{2q+1}(u)}{ u^{2q}}d u
	\end{align*}|
with $\ds C_{m,q}^1=\frac{1}{ 2^{2m}}\sum_{k=0}^m(-1)^{k+q}\binom{2m+1}{ m-k}(2k+1)^{2q-1}$.
But, for $1\leq q\leq m$ the integral $\ds \int_0^\infty\frac{\sin^{2m+1}t}{ t^{2q}}d t$ converges, so that using $iv$ from 
Proposition \ref{pr1}, we conclude that we must have $ C_{m,q}^1=0$, and consequently
$$\int_0^\infty\frac{\sin^{2m+1}t}{ t^{2q}}d t
=\frac{1}{2^{2m}}\sum_{k=0}^m(-1)^{k+q}\binom{2m+1}{ m-k}\frac{(2k+1)^{2q-1}}{(2q-1)!}\ln (2k+1).$$
This concludes the proof of the theorem. \bg
\end{proof}
\begin{corollary}
If $q$ is a positive integer and $k$ is a nonzero integer then 
$\ds\int_0^\infty\frac{\sin^{q+2k}t}{t^q}dt$ is a rational multiple of $\pi$.
\end{corollary}

\begin{example} Here are  some numerical examples for even values of the difference $p-q$:
\[
\def\arraystretch{1.7}
\begin{array}{|c|c|c|c|c|}
\hline
\int_0^\infty\frac{\sin^{p}t}{t^q}dt&p=q&k=q+2&k=q+4&k=q+6\\ \hline
q=1& \frac{\pi }{2} & \frac{\pi }{4} & \frac{3 \pi }{16} & \frac{5 \pi }{32} \\  \hline
q=2& \frac{\pi }{2} & \frac{\pi }{4} & \frac{3 \pi }{16} & \frac{5 \pi }{32} \\  \hline
q=3& \frac{3 \pi }{8} & \frac{5 \pi }{32} & \frac{7 \pi }{64} & \frac{45 \pi }{512} \\  \hline
q=4& \frac{\pi }{3} & \frac{\pi }{8} & \frac{\pi }{12} & \frac{25 \pi }{384} \\  \hline
\end{array}
\]
\end{example}

\begin{example} Here are some numerical examples for odd values of the difference $p-q$, (note that the considered integral diverges for $q=1$ in this case):

\[
\def\arraystretch{1.7}
\begin{array}{|c|c|c|c|c|}
\hline
\int_0^\infty\frac{\sin^{p}t}{t^q}dt&p=q+1&p=q+3&p=q+5\\ \hline
q=2& \frac{3 }{4}\log3 & \frac{15 }{16}\log3-\frac{5}{16}\log5 & \frac{63}{64}\log3-\frac{35}{64}\log5+\frac{7}{64}\log7  \\  \hline
q=3& \log2 & \frac{3 }{2}\log2-\frac{9}{16}\log3 & \frac{9 \pi }{4}\log2-\frac{9}{8}\log3  \\  \hline
q=4& -\frac{45 }{32}\log3+\frac{125}{96}\log5 & -\frac{189 }{128}\log3+\frac{875}{384}\log5-\frac{343}{384}\log7 & 
-\frac{270}{512}\log3+\frac{1500}{512}\log5-\frac{1029}{512}\log7  \\  \hline
q=5& -2\log2+\frac{27 }{16}\log3 &  -5\log2+\frac{27 }{8}\log3 & 
-\frac{55}{6}\log2+\frac{1215}{256}\log3-\frac{625}{768}\log5  \\  \hline
\end{array}
\]
\end{example}
\bg
\begin{remark}
It is interesting to note that these integrals do not appear 
in \cite{grad}, but only particular cases are given in paragraphs 3.82--3.83. 
\end{remark}

\end{document}